# Darboux Slant Ruled Surfaces

## Mehmet Önder, Onur Kaya

*Celal Bayar University, Faculty of Arts and Sciences, Department of Mathematics, Muradiye Campus, 45047 Muradiye, Manisa, Turkey.*
E-mails: mehmet.onder@cbu.edu.tr, onur.kaya@cbu.edu.tr

**Abstract**

In this study, we introduce Darboux slant ruled surfaces in the Euclidean 3-space $E^3$ which is defined by the property that the Darboux vector of orthonormel frame of ruled surface makes a constant angle with a fixed, non-zero direction. We obtain the characterizations of Darboux slant ruled surfaces regarding the conical curvature $\kappa$ and give the relations between Darboux slant ruled surfaces and some other slant ruled surfaces.

## 1. Introduction

Special curves and surfaces which are defined by some properties according to their curvatures are the most fascinating problem of the differential geometry. The well-known types of special curves are helices, involute-evolute curves, Bertrand curves and slant helices. Helix curve is defined by the property that the tangent line of the curve makes a constant angle with a fixed straight line (the axis of the general helix) [2]. The classical result for a helix was stated by Lancret in 1802 and first proved by B. de Saint Venant in 1845: *A necessary and sufficient condition that a curve to be a general helix is that the ratio of the first curvature to the second curvature be constant i.e., $\kappa/\tau$ is constant along the curve, where $\kappa$ and $\tau$ denote the first and second curvatures of the curve, respectively* [11]. Recently, Izumiya and Takeuchi have defined a new special curve called slant helix for which the principal normal lines of the curve make a constant angle with a fixed direction and they have given a characterization of slant helix in the Euclidean 3-space $E^3$ [4]. The spherical images of the Frenet vectors of a slant helix have been studied by Kula and Yaylı and they have obtained that the spherical images of a slant helix are spherical helices [6]. Later, Kula and et al have obtained some new characterizations of slant helices in $E^3$ [7]. Monterde has shown that for a curve with constant curvature and non-constant torsion the principal normal vector makes a constant angle with a fixed constant direction, i.e, the curve is a slant helix [8]. Ali has considered the position vectors of slant helices and obtained some new properties of these curves [1]. Analogue to the definition of slant helix, Önder and et al. have defined $B_2$-slant helix in the Euclidean 4-space $E^4$ by saying that the second binormal vector of a space curve makes a constant angle with a fixed direction and they have given some characterizations of $B_2$-slant helices in the Euclidean 4-space $E^4$ [10]. Then, Gök, Camcı and Hacısalihoğlu have given the generalization of slant helices in the Euclidean $n$-space $E^n$ and called $V_n$-slant helix [3]. Later, Zıplar, Şenol and Yaylı have defined a new type of helices called Darboux helices and given some characterizations for these curves [12].

Moreover, Önder has considered the notion of "slant helix" for ruled surfaces and defined slant ruled surfaces in the Euclidean 3-space $E^3$ by the property that the vectors of the frame of surface along striction curve make constant angles with fixed lines [9]. He has shown that the striction curves of developable slant ruled surfaces are helices or slant helices.

In this work, we introduce Darboux slant ruled surfaces in $E^3$. We give characterizations of these special surfaces and obtain the relationships between Darboux slant ruled surfaces and other slant ruled surfaces.

## 2. Ruled Surfaces in the Euclidean 3-space $E^3$

In this section, we give a brief summary about ruled surfaces in $E^3$.

A ruled surface $S$ is a special surface generated by a continuous movement of a line along a curve and has the parametrization

$$\vec{r}(u,v) = \vec{f}(u) + v\vec{q}(u), \qquad (1)$$

where $\vec{f} = \vec{f}(u)$ is a regular curve in $E^3$ defined on an open interval $I \subset \mathbb{R}$ and $\vec{q} = \vec{q}(u)$ is a unit direction vector of an oriented line in $E^3$. The curve $\vec{f} = \vec{f}(u)$ is called base curve or generating curve of the surface and various positions of the generating lines $\vec{q} = \vec{q}(u)$ are called rulings. In particular, if the direction of $\vec{q}$ is constant, then the ruled surface is said to be cylindrical, and non-cylindrical otherwise.

Let $\vec{m}$ be the unit normal vector of the ruled surface $S$. Then if $v$ decreases infinitely, along a ruling $u = u_1$, the unit normal $\vec{m}$ approaches a limiting direction. This direction is called asymptotic normal (central tangent) direction and is defined by

$$\vec{a} = \lim_{v \to \pm\infty} \vec{m}(u_1, v) = \frac{\vec{q} \times \dot{\vec{q}}}{\|\dot{\vec{q}}\|}. \qquad (2)$$

The point at which the unit normal of $S$ is perpendicular to $\vec{a}$ is called the striction point (or central point) $C$ and the set of striction points of all rulings is called striction curve of the surface. The parametrization of the striction curve $\vec{c} = \vec{c}(u)$ on a ruled surface is given by

$$\vec{c}(u) = \vec{f} - \frac{\langle \dot{\vec{q}}, \dot{\vec{f}} \rangle}{\langle \dot{\vec{q}}, \dot{\vec{q}} \rangle} \vec{q}. \qquad (3)$$

The vector $\vec{h}$ defined by $\vec{h} = \vec{a} \times \vec{q}$ is called central normal vector. Then the orthonormal system $\{C; \vec{q}, \vec{h}, \vec{a}\}$ is called Frenet frame of the ruled surface $S$ where $C$ is the central point and $\vec{q}, \vec{h}, \vec{a}$ are unit vectors of ruling, central normal and central tangent, respectively.

The set of all bound vectors $\vec{q}(u)$ at the point $O$ constitutes a cone which is called *directing cone* of the ruled surface $S$. The end points of unit vectors $\vec{q}(u)$ drive a spherical curve $k_1$ on the unit sphere $S_1^2$ and this curve is called *spherical image* of ruled surface $S$, whose arc length is denoted by $s_1$.

For the Frenet formulae of the ruled surface $S$ and of its directing cone with respect to the arc length $s_1$ we have

$$\begin{bmatrix} d\vec{q}/ds_1 \\ d\vec{h}/ds_1 \\ d\vec{a}/ds_1 \end{bmatrix} = \begin{bmatrix} 0 & 1 & 0 \\ -1 & 0 & \kappa \\ 0 & -\kappa & 0 \end{bmatrix} \begin{bmatrix} \vec{q} \\ \vec{h} \\ \vec{a} \end{bmatrix}, \qquad (4)$$

where $\kappa$ is called conical curvature of the directing cone (For details [5]). The Frenet formulae can be interpreted kinematically as follows: If $\vec{q}$ traverses the directing cone in such a way that $s$ is the time parameter, then the moving frame $\{C; \vec{q}, \vec{h}, \vec{a}\}$ moves in accordance with (5). This motion contains, apart from an instantaneous translation, and instantaneous rotation with angular velocity vector given by the Darboux vector

$$\vec{W} = \kappa \vec{q} + \vec{a}. \qquad (5)$$

The direction of the Darboux vector is that of instantaneous axis of rotation, and its length $\|\vec{W}\| = \sqrt{1 + \kappa^2}$ is the scalar angular velocity. Then, Frenet formulae (5) can be given as follows,

$$\vec{q}' = \vec{W} \times \vec{q}, \qquad \vec{h}' = \vec{W} \times \vec{h}, \qquad \vec{a}' = \vec{W} \times \vec{a}. \qquad (6)$$

**Definition 2.1.** ([9]) Let $S$ be a regular ruled surface in $E^3$ given by the parametrization

$$\vec{r}(s,v) = \vec{c}(s) + v\vec{q}(s), \quad \|\vec{q}(s)\| = 1,$$

where $\vec{c}(s)$ is striction curve of $S$ and $s$ is arc length parameter of $\vec{c}(s)$. Let the Frenet frame of $S$ be $\{\vec{q}, \vec{h}, \vec{a}\}$. Then $S$ is called a $q$-slant ($h$-slant or $a$-slant, respectively) ruled surface if the ruling (the vector $\vec{h}$ or the vector $\vec{a}$, respectively) makes a constant angle with a fixed non-zero direction $\vec{u}$ in the space, i.e.,

$$\langle \vec{q}, \vec{u} \rangle = \cos\theta = \text{constant}; \ \theta \neq \frac{\pi}{2},$$

$$(\langle \vec{h}, \vec{u} \rangle = \cos\theta = \text{constant}; \ \theta \neq \frac{\pi}{2} \text{ or } \langle \vec{a}, \vec{u} \rangle = \cos\theta = \text{constant}; \ \theta \neq \frac{\pi}{2}, \text{ respectively}).$$

In [9], Önder has given the characterizations of slant ruled surfaces according to Frenet frame of a ruled surface which is defined by arc length parameter $s$ of striction curve. In this paper, we will study the Darboux slant ruled surfaces by considering Eq. (5) and give the relationships between slant ruled surfaces and Darboux slant ruled surfaces. So, first we prove the following theorem for $\vec{h}$-slant ruled surfaces.

**Theorem 2.1**. *Let $S$ be a regular ruled surface in $E^3$ with Frenet frame $\{\vec{q}, \vec{h}, \vec{a}\}$ and conical curvature $\kappa \neq 0$. Then $S$ is an $h$-slant ruled surface if and only if the function*

$$\frac{\kappa'}{\left(1+\kappa^2\right)^{3/2}}, \tag{7}$$

*is constant.*

***Proof:*** Assume that $S$ is an $h$-slant ruled surface in $E^3$. So, for a non-zero constant $c \in \mathbb{R}$ we can write

$$\langle \vec{h}, \vec{u} \rangle = c.$$

where $\vec{u}$ is a non-zero, fixed direction. Then, for the vector $\vec{u}$ we have

$$\vec{u} = b_1(s_1)\vec{q}(s_1) + c\vec{h}(s_1) + b_2(s_1)\vec{a}(s_1) \tag{8}$$

where $b_1 = b_1(s_1)$ and $b_2 = b_2(s_1)$ are smooth functions of arc length parameter $s_1$. On the other hand, $\vec{u}$ is a fixed direction, that is $\vec{u}' = 0$. Since the Frenet frame $\{\vec{q}, \vec{h}, \vec{a}\}$ is linearly independent, differentiation of (8) gives

$$\begin{cases} b_1' - c = 0, \\ b_1 - \kappa b_2 = 0, \\ b_2' + c\kappa = 0. \end{cases} \tag{9}$$

From the second equation of system (9) we have

$$b_1 = \kappa b_2. \tag{10}$$

Moreover, since $\|\vec{u}\|$ is constant along striction curve, we obtain

$$b_1^2 + c^2 + b_2^2 = constant. \tag{11}$$

Substituting (10) in (11) gives

$$b_2^2 \left(1 + \kappa^2\right) = n^2 = constant. \tag{12}$$

If $n = 0$, then $b_2 = 0$ and from (9) we have $b_1 = 0$, $c = 0$. This means that $\vec{u} = 0$ which is a contradiction. Thus, $n \neq 0$. Then from (12) it is obtained that

$$b_2 = \pm \frac{n}{\sqrt{1+\kappa^2}}. \tag{13}$$

Considering the third equation of the system (9), from (13) we have

$$\frac{d}{ds_1}\left[\pm \frac{n}{\sqrt{1+\kappa^2}}\right] = -c\kappa. \tag{14}$$

For a non-zero real constant $d$, this can be written as

$$\frac{\kappa'}{(1+\kappa^2)^{3/2}} = \frac{c}{n} = d,$$

which is desired.

Conversely, assume that the function in (7) is constant, i.e.,

$$\frac{\kappa'}{(1+\kappa^2)^{3/2}} = d = constant.$$

We define the function

$$\vec{u} = \frac{\kappa}{\sqrt{1+\kappa^2}}\vec{q} + d\vec{h} + \frac{1}{\sqrt{1+\kappa^2}}\vec{a}. \tag{15}$$

From (15) and Frenet formulae (4), we have $\vec{u}' = 0$, i.e, $\vec{u}$ is a constant vector. On the other hand,

$$\langle \vec{h}, \vec{u} \rangle = d = constant.$$

Thus, $S$ is an $h$-slant ruled surface in $E^3$.

## 3. Darboux Slant Ruled Surfaces in the Euclidean 3-space $E^3$

In this section, we consider the notion of "slant" for Darboux vector $\vec{W} = \kappa\vec{q} + \vec{a}$ and give some theorems for Darboux slant ruled surfaces in the Euclidean 3-space. First, we give the following definition.

**Definition 3.1.** Let $S$ be a regular ruled surface in $E^3$ given by the parametrization

$$\vec{r}(s,v) = \vec{c}(s) + v\vec{q}(s), \quad \|\vec{q}(s)\| = 1, \tag{16}$$

where $\vec{c}(s)$ is the striction curve of $S$ and $s$ is the arc length parameter of $\vec{c}(s)$. Let the orthonormal frame $\{\vec{q}, \vec{h}, \vec{a}\}$ and the function $\kappa = \kappa(s)$ denote the Frenet frame and conical curvature of $S$, respectively. Then, $S$ is called a Darboux slant ruled surface if its Darboux vector $\vec{W}$ makes a constant angle $\varphi$ with a non-zero, fixed direction $\vec{u}$ in the space, i.e.,

$$\langle \vec{W}, \vec{u} \rangle = \cos\varphi.$$

Then, we give the following theorems for Darboux slant ruled surfaces. Whenever we talk about $S$, we will mean that the surface is regular and has the parametrization and Frenet elements as assumed in Definition 3.1.

**Theorem 3.1.** *Let $S$ be a ruled surface in $E^3$. If $S$ is a Darboux slant ruled surface, then the conical curvature $\kappa$ is constant.*

***Proof:*** Let $S$ be a Darboux slant ruled surface. Then for a non-zero, fixed direction $\vec{u}$, we have,

$$\langle W, \vec{u} \rangle = \text{constant} \tag{17}$$

Taking the derivative of (17) gives

$$\langle W', \vec{u} \rangle = 0. \tag{18}$$

From (18) and Frenet formulae (4), we can write,

$$\kappa' \langle \vec{q}, d \rangle = 0 \tag{19}$$

From (19) we get two possibilities:

$$\begin{cases} \kappa = \text{constant}, \\ \langle \vec{q}, \vec{u} \rangle = 0. \end{cases} \tag{20}$$

If $\langle \vec{q}, \vec{u} \rangle = 0$, then $\vec{u}$ is perpendicular to the vector $\vec{q}$ and can be written as

$$\vec{u} = a_1 \vec{h} + a_2 \vec{a}. \tag{21}$$

where $a_1$, $a_2$ are smooth functions of $s_1$. By taking the derivative of (21) and using the Frenet formulae given in (4) we have

$$-a_1 \vec{q} + (a_1' - \kappa a_2)\vec{h} + (a_2' + \kappa a_1)\vec{a} = 0. \tag{22}$$

Considering that the Frenet frame $\{\vec{q}, \vec{h}, \vec{a}\}$ is linearly independent, from (22) we obtain the following system

$$\begin{cases} a_1 = 0, \\ a_1' - \kappa a_2 = 0, \\ a_2' + \kappa a_1 = 0. \end{cases} \tag{23}$$

From (23) we get $a_1 = a_2 = 0$ which gives us that $\vec{u} = 0$ which is a contradiction, that is $\langle \vec{q}, \vec{u} \rangle \neq 0$. Therefore, $\kappa = \text{constant}$.

***Corollary 3.1.*** Let $S$ be a regular ruled surface in $E^3$. If $S$ is a Darboux slant ruled surface then $\det(\vec{W}, \vec{W}', \vec{W}'') = 0$.

***Proof:*** From Darboux vector and its derivatives we have,

$$\vec{W} = \kappa \vec{q} + \vec{a},$$
$$\vec{W}' = \kappa' \vec{q},$$
$$\vec{W}'' = \kappa'' \vec{q} + \kappa' \vec{h},$$

and so the determinant would be,

$$\det(\vec{W}, \vec{W}', \vec{W}'') = (\kappa')^2.$$

If $S$ is a Darboux slant ruled surface, from Theorem 3.1 we have that $\kappa =$ constant which gives us $\det(\vec{W}, \vec{W}', \vec{W}'') = 0$.

**Theorem 3.2.** *Every $h$-slant ruled surface is also a Darboux slant ruled surface.*

**Proof:** Let $S$ be an $h$-slant ruled surface. So, we have

$$\langle \vec{h}, \vec{u} \rangle = \cos\theta = \text{constant}, \tag{24}$$

where $\theta$ is constant angle between the vectors $\vec{h}, \vec{u}$. From Theorem 2.1, we know that the axis of $h$-slant ruled surface is given by

$$\vec{u} = \frac{\kappa}{\sqrt{1+\kappa^2}} \vec{q} + \cos\theta \vec{h} + \frac{1}{\sqrt{1+\kappa^2}} \vec{a}. \tag{25}$$

By using (25) and considering Darboux vector $\vec{W} = \kappa \vec{q} + \vec{a}$ we can write,

$$\vec{u} = \frac{\vec{W}}{\|\vec{W}\|} + \cos\theta \vec{h}, \tag{26}$$

which gives that $\|\vec{u}\| = \sqrt{1 + \cos^2\theta}$ is constant. Then, from (26) we obtain

$$\langle \vec{W}, \vec{u} \rangle = \|\vec{W}\| + \cos\theta \langle \vec{W}, \vec{h} \rangle. \tag{27}$$

Furthermore, from vector calculus we know that $\langle \vec{W}, \vec{u} \rangle = \|\vec{W}\| \|\vec{u}\| \cos\lambda$ where $\cos\lambda$ is the angle between $\vec{W}$ and $\vec{u}$. Then from (27) it follows

$$\|\vec{u}\| \cos\lambda = 1. \tag{28}$$

So, from (28) we obtain,

$$\cos\lambda = \frac{1}{\sqrt{1 + \cos^2\theta}} = \text{constant},$$

which means that $\langle \vec{W}, \vec{u} \rangle = $ constant. Thus, the surface $S$ is a Darboux slant ruled surface.

***Theorem 3.3.*** Let $S$ be a Darboux slant ruled surface. If $S$ is an $h$-slant ruled surface then $S$ is also a $q$-slant (or $a$-slant) ruled surface.

***Proof:*** Since $S$ is a Darboux slant ruled surface, from Theorem 3.1 we have that the conical curvature $\kappa$ is constant. If $S$ is an $h$-slant ruled surface, then from (15) the axis of the surface is

$$\vec{u} = \frac{\kappa}{\sqrt{1+\kappa^2}}\vec{q} + d\vec{h} + \frac{1}{\sqrt{1+\kappa^2}}\vec{a}.$$

Since $\kappa$ is constant, the last equality gives that $\langle \vec{q}, \vec{u} \rangle =$ constant, (or $\langle \vec{a}, \vec{u} \rangle =$ constant), i.e., $S$ is also a $q$-slant (or $a$-slant) ruled surface.

The converse of Theorem 3.3 is not always true. It is valid for a special value of coefficient $a_3$. So, we give the following special case.

***Theorem 3.4.*** Let $S$ be a Darboux slant ruled surface. Then, $S$ is an $h$-slant ruled surface if and only if $S$ is also an $a$-slant ruled surface with constant angle $\lambda$ defined by $\langle \vec{a}, \vec{u} \rangle = \cos \lambda = \frac{\cos \varphi}{1+\kappa^2}$ where $\vec{u}$ is a non-zero, fixed vector and $\cos \varphi = \langle \vec{W}, \vec{u} \rangle$.

***Proof:*** Since $S$ is a Darboux slant ruled surface, we have

$$\langle \vec{W}, \vec{u} \rangle = \cos \varphi = \text{constant}. \tag{29}$$

where $\vec{u}$ is a non-zero, fixed vector and $\varphi$ is the constant angle between the vectors $\vec{W}, \vec{u}$. For the vector $\vec{u}$ we can write,

$$\vec{u} = a_1 \vec{q} + a_2 \vec{h} + a_3 \vec{a}. \tag{30}$$

where $a_i = a_i(s_1)$, $(1 \leq i \leq 3)$ are smooth functions of $s_1$. By taking the derivative of (30) we obtain,

$$(a_1' - a_2)\vec{q} + (a_1 + a_2' - a_3\kappa)\vec{h} + (a_2\kappa + a_3')\vec{a} = 0. \tag{31}$$

Since the Frenet frame $\{\vec{q}, \vec{h}, \vec{a}\}$ is linearly independent, we have the following system,

$$\begin{cases} a_1' - a_2 = 0, \\ a_1 + a_2' - a_3\kappa = 0, \\ a_2\kappa + a_3' = 0. \end{cases} \tag{32}$$

Writing (30) in (29) gives

$$\kappa a_1 + a_3 = \cos\varphi. \tag{33}$$

By using (33) and the second equation of (32) we have,

$$a_2' = \frac{(1+\kappa^2)a_3 - \cos\varphi}{\kappa}. \tag{34}$$

Since $S$ is a Darboux slant ruled surface, from Theorem 3.1 we have that $\kappa$ is constant. Now, from (34) we have that $a_2$ is constant if and only if $a_3$ is a constant given by $a_3 = \frac{\cos\varphi}{1+\kappa^2}$, i.e, $S$ is an $h$-slant ruled surface if and only if $S$ is also an $a$-slant ruled surface with constant angle $\lambda$ defined by $\langle \vec{a}, \vec{u} \rangle = \cos\lambda = \frac{\cos\varphi}{1+\kappa^2}$.